\documentclass[a4paper,11pt]{amsart}
\usepackage{amsmath}
\DeclareMathOperator{\sw}{s}
\DeclareMathOperator{\invol}{{\mathfrak{Invol}}}
\DeclareMathOperator{\Ad}{Ad}
\DeclareMathOperator{\id}{id}
\DeclareMathOperator{\inv}{inv}
\DeclareMathOperator{\sgn}{sgn}
\DeclareMathOperator{\Ind}{Ind}
\DeclareMathOperator{\idem}{idem}
\newcommand{\ootimes}{{\otimes}}

\newcommand{\simalg}{\simeq_{\mathrm{alg}}}
\newcommand{\marginextend}[1]{ \addtolength{\oddsidemargin}{-#1}  \addtolength{\evensidemargin}{-#1}
  \addtolength{\textwidth}{#1}\addtolength{\textwidth}{#1}}
\newcommand{\updownextend}[1]{ \addtolength{\topmargin}{-#1}  \addtolength{\textheight}{#1}
\addtolength{\textheight}{#1}}
\marginextend{12.5mm}
\updownextend{0.0mm}
\theoremstyle{definition}
\newtheorem{point}{}[section]

\newtheorem{remark}[point]{Remark}
\theoremstyle{plain}
\newtheorem{cor}[point]{Corollary}
\newtheorem{theorem}[point]{Theorem}
\newtheorem{lemma}[point]{Lemma}
\begin{document}
\title{On the $H$-ring structure of infinite Grassmannians}
\author{Gyula Lakos}
\address{Institute of Mathematics, E\"otv\"os University, P\'azm\'any P\'eter s.~1/C,  Budapest, H--1117, Hungary}
\email{lakos@cs.elte.hu}
\thanks{The author thanks for the hospitality of the Alfr\'ed R\'enyi Institute of Mathematics,
where the idea of this paper was conceived.}
\begin{abstract}
The $H$-ring structure of certain infinite dimensional Grassmannians is discussed using various
algebraic and analytical methods but avoiding cellular arguments.
These methods allow us to treat these Grassmannians in greater generality.
\end{abstract}
\maketitle
\section*{Introduction}
Infinite dimensional Grassmannians are often used as realizations of classifying spaces like
$\mathbb Z\times BO$ and $\mathbb Z\times BU$.
The two spaces mentioned above classify $KO_0$ and $K_0$ respectively; consequently they possess $H$-ring
structures induced from the direct sum, switch, and tensor product constructions of virtual vector
bundles.
In standard textbook constructions, the homotopy additive structure of these Grassmannians is often described
in explicit terms, see e.~g.~\cite{SW}.
However, the existence of the homotopy product map is rather inferred from
principles of homotopy theory instead of constructed explicitly.
The main tools are approximation, weak equivalences, and universality.
One might get the impression that cellular arguments are unavoidable in that respect.
But the truth is that the $H$-ring structure of the infinite Grassmannians is much more an algebraic,
or perhaps analytic, matter than a combinatorial one.
Our objective here is to work out this structure without obscuring it with  cellular topology.
First, we discuss the algebraic $H$-ring structure.
We take a ring $\mathfrak A$ endowed by a polymetric structure.
Then we consider the virtual Grassmannian $\mathcal G^{(2)}(\mathfrak A)$ and
the ordinary infinite Grassmannian $\mathcal G(\mathfrak A)$ associated to $\mathfrak A$.
We will show that they possess commutative unital $H$-ring structures.
Strictly speaking, this holds if $\mathfrak A$ is a commutative ring, but one
can formulate this phenomenon in terms of tensor products such that it applies more generally.
Second, if $\mathfrak A$ is a locally convex algebra, then the algebraic $H$-ring structure
implies a topological, in fact, smooth, $H$-ring structure.
Third, if $\mathfrak A$ satisfies somewhat stronger conditions, then the smooth
$H$-structure implies a smooth algebraic $H$-structure without stabilization.
This all applies to $\mathfrak A=\mathbb R$ and $\mathfrak A=\mathbb C$, corresponding
to the classical cases mentioned above. We conclude the paper with a
notion of dimension, and the discussion of Fredholm operators.

\section{Polymetric rings and associated matrix spaces}

We say that the set $\Omega$ is an infinite  set of polynomial growth (spg)
if it is endowed by a set of real valued  functions $\mathcal S^{-\infty}_\Omega$ on $\Omega$
such that there is a bijection $q:\Omega\rightarrow \mathbb N$ so that
$\mathcal S^{-\infty}_\Omega=q^*\mathcal S^{-\infty}_{\mathbb N}$, where $\mathcal S^{-\infty}_{\mathbb N}$
is the set of real valued functions of at most polynomial growth on $\mathbb N$.
We say that $\Omega$ is a finite set of polynomial growth
if it is finite, and it is endowed with the set of arbitrary real valued functions $\mathcal S^{-\infty}_\Omega$.
If $\Omega_1$ and $\Omega_2$ are spg's, then one can naturally construct the spg's
$\Omega_1\,\dot\cup\,\Omega_2$ and $\Omega_1\times\Omega_1$.
If $\Omega=\Omega_1\,\dot\cup\,\Omega_2$ as spg's, then
we say that $\Omega$ decomposes to $\Omega_1$ and $\Omega_2$.
(An arbitrary set-theoretical decomposition is not sufficient in general.)

We say that the topological ring $\mathfrak A$ is a polymetric ring if
\begin{itemize}
\item[a.)] its topology is induced by a family of ``seminorms'' $p:\mathfrak A\rightarrow[0,+\infty)$ such that
$p(0)=0$, $p(-X)=p(X)$, $p(X+Y)\leq p(X)+p(Y)$;

\item[b.)]for each ``seminorm'' $p$ there exists a ``seminorm'' $\tilde p$ such that
$p(XY)\leq\tilde p(X)\tilde p(X)$ holds.
\end{itemize}

This is a large class of topological rings: it includes locally convex algebras just as
discrete rings with the convention $p(X)=1$ for $X\neq 0$.
In what follows, $\mathfrak A$ is assumed to be a  separated, sequentially complete polymetric ring.

Suppose that $\mathfrak A$ is a polymetric ring, $\Omega$ is an spg.
Then we can define the algebra of rapidly decreasing matrices $\mathcal K_{\Omega}(\mathfrak A)$ and
the algebra of matrices of pseudodifferential size $\Psi_{\Omega}(\mathfrak A)$; essentially as in \cite{SS}.
(In the special case when $\mathfrak A$ is the discrete ring,  the space
$\mathcal K_{\Omega}(\mathfrak A)$ is the space of matrices with finitely many non-zero
elements, and $\Psi_{\Omega}(\mathfrak A)$ is the space of matrices such that
every column and row has only finitely many non-zero elements.)
More generally, we can take  the spaces of $\Omega'\pmb\times\Omega$ matrices
$\Psi_{\Omega',\Omega}(\mathfrak A)$ and $\mathcal K_{\Omega',\Omega}(\mathfrak A)$.
(This bold $\pmb\times$ is reserved for matrix shape.)
If $\mathfrak A$ is a polymetric ring, then we may take its unital extension $\mathfrak A^+$.
We will be a bit vague about this construction:
In the general case, it may be
the group $\mathbb Z\oplus \mathfrak A$ endowed with the naturally extended structure,
but if $\mathfrak A$ is a locally convex algebra over $\mathbb K$, then there is
no danger in taking $\mathbb K\oplus\mathfrak A$ with the naturally extended algebra structure.
Let, in general, $1_\Omega=\sum_{\omega\in\Omega}\mathbf e_{\omega,\omega}\in \Psi_{\Omega}(\mathfrak A^+)$;
and let, in general, $0_\Omega$ denote the nullmatrix over $\Omega$.
The matrix $A\in \Psi_{\Omega',\Omega}(\mathfrak A^+)$ is invertible if there is an element
$B\in \Psi_{\Omega,\Omega'}(\mathfrak A^+)$ such that $AB=1_{\Omega'}$ and $BA=1_{\Omega}$.
We denote the set of those by $\Psi_{\Omega',\Omega}(\mathfrak A^+)^\star$, and call the them units.
The unit group $\mathcal K_{\Omega}(\mathfrak A)^\star$ is the group
of invertible elements $1_{\Omega}+A\in \Psi_{\Omega}(\mathfrak A^+)$,
where $A\in \mathcal K_{\Omega}(\mathfrak A)$,
but with topology induced from $\mathcal K_{\Omega}(\mathfrak A)$.
Furthermore, let $\Psi^{(2)}_{\Omega',\Omega}(\mathfrak A^+)$
be the space of pairs $\langle B,A\rangle$, where $A,B\in \Psi_{\Omega',\Omega}(\mathfrak A^+)$,
but $B-A\in \mathcal K_{\Omega',\Omega}(\mathfrak A)$.
Its topology is  induced jointly from $A,B$ and $B-A$ with respect to the appropriate spaces, respectively.
Then the unit group $\Psi^{(2)}_{\Omega',\Omega}(\mathfrak A^+)^\star$ can be taken.
Conjugation induces a continuous map
$\Ad: \Psi_{\Omega',\Omega}(\mathfrak A^+)^\star\times \mathcal K_{\Omega}(\mathfrak A)^\star
\rightarrow \mathcal K_{\Omega'}(\mathfrak A)^\star$, etc.

Some simple invertible matrices in $\Psi_{\Omega',\Omega}(\mathfrak A^+)$ are
as follows.
Let $r:\Omega\rightarrow \Omega'$ be an isomorphism of spg's.
Then we take
$\hat r=\sum_{\omega\in\Omega} \mathbf e_{r(\omega),\omega}\in\Psi_{\Omega',\Omega}(\mathfrak A^+)^\star$.
For $A\in \Psi_{\Omega}(\mathfrak A^+)$ the map
$r_*: A \mapsto \hat rA\hat r^\top$ has the effect that
$r_*\left(\sum_{n,m\in\Omega} a_{n,m}\mathbf e_{n,m}\right)=\sum_{n,m\in\Omega} a_{n,m}\mathbf e_{r(n),r(m)}$.
Hence we call such $r_*$ isomorphic relabeling maps.
More generally, if $r:\Omega\rightarrow \Omega'$ is only
a map of spg's such that $r(\Omega)$ and $\Omega'\setminus r(\Omega)$ decomposes
$\Omega'$ as spg's, and $r$ induces an isomorphism between the spg structure of
$\Omega$ and the one of $\Omega$ restricted to $r(\Omega)$,
then $\hat r$ and $r_*$ can be taken. We still call these $r_*$ (not necessarily isomorphic) relabeling maps.
Another  natural operation is the direct sum of matrices. For example,
if $A\in  \Psi_{\Omega}(\mathfrak A), B\in  \Psi_{\Omega'}(\mathfrak A)$, then we can consider
$A\oplus B\in  \Psi_{\Omega\dot\cup\Omega'}(\mathfrak A)$ which is a colloquial notation for the block matrix
$\begin{bmatrix} A&\\&B  \end{bmatrix}\in  \Psi_{\Omega\dot\cup\Omega'}(\mathfrak A)$.
Indexing direct sums might be confusing, especially if the construction is iterated.
We take disjoint union for the index set, but if we take
direct sum of matrices with the same index set $\Omega$, then we might use
$\Omega$, $\Omega'$, $\Omega''$, etc. for the components, or
$\{0\}\times\Omega$, $\{1\}\times\Omega$, $\{2\}\times\Omega$, etc., depending on the situation.
Later, when we consider block matrices of (matrices indexed by $\Omega$) indexed by $\Gamma$
then we consider those matrices as matrices indexed by $\Gamma\times \Omega$.
A $\Gamma_2\pmb\times \Gamma_1$ block matrix of $\Omega_2\pmb\times \Omega_1$
matrices will be considered as a $\Gamma_2\times \Omega_2\pmb\times \Gamma_1\times\Omega_1$
matrix.

We apply the following notational conventions in unital polymetric rings:

(a) We write $\bar a$ for $1-a$. (For $a\in\Psi_{\Omega}(\mathfrak A^+)$ it is, of course, $\bar a=1_\Omega-a$.)

(b) For  $\Xi\times \Xi$ matrices over $\Psi_{\Omega}(\mathfrak A^+)$, we use the notation
\[\sw_{n,m}(a)=1_{(\Xi\setminus\{n,m\})\times\Omega}
+\bar a\mathbf e_{n,n}+a\mathbf e_{n,m}+a\mathbf e_{m,n}+\bar a\mathbf e_{m,m}.\]
We see that $\sw_{n,m}(a)$ is an involution if $a$ is an idempotent.
It means that we have a partial switch between  the $n,m$ positions.

(c) Sometimes we write $A^g$ for  $gAg^{-1}$.
We use the abbreviation $A\xrightarrow{g} B$ for $gAg^{-1}=B$.
In fact, sometimes we say that $g$ is a morphism between $A$ and $B$.

(d) If $b-a\in \mathcal K_\Omega(\mathfrak A)$, then we use the notation $a\approx b$.

(e) $\invol(\mathfrak A)$ denotes the subspace of involutions in $\mathfrak A$.

\section{Virtual Grassmannians and natural operations on them}
\begin{point}
We  define the virtual Grassmannian $\mathcal G^{(2)}_{\Omega}(\mathfrak A)$
as  $\invol(\Psi^{(2)}_{\Omega}(\mathfrak A^+))$,
i.~e. the space of pairs $\langle b,a\rangle$ such that $a,b\in \Psi_{\Omega}(\mathfrak A^+)$,
$b-a\in \mathcal K_{\Omega}(\mathfrak A)$, and $a,b$ are idempotents.
For us, such  pairs are  virtual idempotents $b$``$-$''$a$;
we will call them as pairs of idempotents.
We refer to the first term as the leading term, and we refer to the second term as the base term.

(a) The sum of pairs operation is defined for $\langle b,a \rangle\in \mathcal G^{(2)}_{\Omega}(\mathfrak A)$
and $\langle d,c \rangle\in \mathcal G^{(2)}_{\Xi}(\mathfrak A)$ as
\[\langle b,a\rangle\oplus \langle d,c\rangle:=\langle b\oplus d,a\oplus c\rangle
\in \mathcal G^{(2)}_{\Omega\dot\cup\Xi}(\mathfrak A).\]

(b) We define the inverse pair for $\langle b,a \rangle\in \mathcal G^{(2)}_{\Omega}(\mathfrak A)$ as
\[\langle b,a\rangle^{\mathrm{inv}}:=\langle\bar b,\bar a\rangle\in\mathcal G^{(2)}_{\Omega}(\mathfrak A) .\]

(c) A very special element is the pair $\mathbf 0:=\langle*,*\rangle$ of $0\times 0$, i.~e. empty, matrices.
We call this as the additive neutral element.

These operations satisfy the natural additive associative, commutative, and neutral element identities:
\[(\langle b,a\rangle\oplus \langle d,c\rangle)\oplus \langle f,e\rangle
\simeq\langle b,a\rangle\oplus (\langle d,c\rangle\oplus \langle f,e\rangle),\]
\[\langle b,a\rangle\oplus \langle d,c\rangle\simeq \langle d,c\rangle\oplus \langle b,a\rangle,\]
\[\langle b,a\rangle \simeq \langle b,a\rangle\oplus\mathbf 0\simeq\mathbf 0\oplus \langle b,a\rangle,\]
where ``$\simeq$'' means that we have equality after we make natural identifications in the index sets,
i.~e.~after particularly simple isomorphic relabelings.
What is apparently lacked is a natural additive inverse element identity.
\end{point}
\begin{point}
Let  $\otimes$ be a suitable tensor product operation of rings.
Again, we will be somewhat vague about the meaning of this term:
In general, we may mean a projective tensor product of polymetric rings,
but, if $\mathfrak A$ is commutative, then we may also consider the ordinary product as
a tensor product operation, i.~e. tensor product over itself.

(d)  For $\langle b,a \rangle\in \mathcal G^{(2)}_{\Omega}(\mathfrak A)$,
$\langle d,c \rangle\in \mathcal G^{(2)}_{\Xi}(\mathfrak B)$,
we define the products of pairs of involutions as
\[\langle b,a\rangle\,\overleftarrow{\otimes}\,\langle d,c\rangle:=
\langle b\otimes d +\bar b\otimes c, a\otimes d +\bar a\otimes c  \rangle
\in \mathcal G^{(2)}_{\Omega\times\Xi}(\mathfrak A\otimes \mathfrak B),\]
and
\[\langle b,a\rangle\,\overrightarrow{\otimes}\,\langle d,c\rangle:=
\langle b\otimes d +a\otimes \bar d, b\otimes c +a\otimes \bar c\rangle
\in \mathcal G^{(2)}_{\Omega\times\Xi}(\mathfrak A\otimes \mathfrak B).\]
So, we have two natural product operations, which may be somewhat strange.

(e) Another special element is $\mathbf 1:=\langle1,0\rangle$, where the elements are $1\times 1$ matrices,
or rather ``scalars''. Again, we will be vague about the ring it is over, we may mean $\mathbb Z$
or the base field $\mathbb K$ of an algebra.

The natural multiplicative associativity, distributive, and neutral element rules hold:
\[(\langle b,a\rangle\,\overleftarrow{\otimes}\,\langle d,c\rangle)\,\overleftarrow{\otimes}
\,\langle f,e\rangle\simeq
\langle b,a\rangle\,\overleftarrow{\otimes}\,(\langle d,c\rangle\,\overleftarrow{\otimes}\,\langle f,e\rangle),\]
\[\langle b,a\rangle\,\overleftarrow{\otimes}\,(\langle d,c\rangle\oplus \langle f,e\rangle)
\simeq (\langle b,a\rangle\,\overleftarrow{\otimes}\,\langle d,c\rangle)\oplus
(\langle b,a\rangle\,\overleftarrow{\otimes}\,\langle f,e\rangle),\]
\[(\langle b,a\rangle\oplus\langle d,c\rangle) \,\overleftarrow{\otimes}\,\langle f,e\rangle)
\simeq (\langle b,a\rangle\,\overleftarrow{\otimes}\,\langle f,e\rangle)\oplus
(\langle d,c\rangle\,\overleftarrow{\otimes}\,\langle f,e\rangle),\]
\[\langle b,a\rangle \simeq \langle b,a\rangle\,\overleftarrow{\otimes}\,\mathbf 1
\simeq\mathbf 1\,\overleftarrow{\otimes}\, \langle b,a\rangle;\]
and similarly for the other product.
We also see that
\[\langle b,a\rangle\,\overleftarrow{\otimes}\,\langle d,c\rangle
\simeq \langle d,c\rangle\,\overrightarrow{\otimes}\,\langle b,a\rangle.\]
What we clearly miss is the equivalence of the two product operations, which is  equivalent to the
problem of multiplicative commutativity.
\end{point}
\begin{point}
We must look for a weaker equivalence relation in order to get the missing identities.
That will be a notion of algebraic homotopy.
Let us define $\mathbf 0_\Xi=\langle0_\Xi,0_\Xi\rangle$
and $\mathbf 0_\Xi'=\langle1_\Xi,1_\Xi\rangle$, in general.
We will also use the notation $\mathbf I_\Xi=\langle1_\Xi,1_\Xi\rangle$ but in different context.
We say that the maps $f_1:X\rightarrow \mathcal G^{(2)}_{\Omega_1}(\mathfrak A)$ and
$f_2:X\rightarrow \mathcal G^{(2)}_{\Omega_2}(\mathfrak A)$ are algebraically homotopic
if there are index sets $\Omega_{10},\Omega_{11},\Omega_{\tilde 20},\Omega_{\tilde 21},$ and a map
\[g:X\rightarrow \Psi^{(2)}_{\Omega_2\dot\cup\Omega_{\tilde 20}\dot\cup\Omega_{\tilde 21},
\Omega_1\dot\cup\Omega_{10}\dot\cup\Omega_{11}}
(\mathfrak A)\]
which allows a multiplicative inverse $g^{-1}$ such that
\[f_1\oplus \mathbf 0_{\Omega_{10}}\oplus\mathbf 0_{\Omega_{11}}'
\xrightarrow{g}f_2\oplus \mathbf 0_{\Omega_{\tilde 20}}\oplus\mathbf 0_{\Omega_{\tilde 21}}';\]
i.~e. if after stabilization the values are conjugate.
We denote this as $f_1\simalg f_2$. This is an equivalence relation:
if $f_2\simalg f_3$ is realized by a similar map $h$, then
\[f_1\oplus \mathbf 0_{\Omega_{10}\dot\cup\Omega_{20}}\oplus\mathbf 0_{\Omega_{11}\dot\cup\Omega_{21}}'
\xrightarrow{(h\oplus\mathbf I_{\Omega_{\tilde 20}\dot\cup\Omega_{\tilde 21}})
(g\oplus\mathbf I_{\Omega_{20}\dot\cup\Omega_{21}})}
f_3\oplus \mathbf 0_{\Omega_{\tilde 20}\dot\cup
\Omega_{\tilde 30}}\oplus\mathbf 0_{\Omega_{\tilde 21}\dot\cup\Omega_{\tilde 31}}'\]
indeed. So, algebraic homotopy is a combination of equivalence by stabilization
and equivalence by conjugation.
In what follows, dependence on $X$ will often be suppressed.

Using this notion of equivalence, which generalizes natural equivalence, we may set out to
demonstrate the missing identities of additive inverse and equality of products.
However, weakening the equivalence relation introduces further problems.
When we consider the generalization of classical structures up to homotopy,
we  have to worry  not only about the classical identities
\[\mathrm{Expr}_1(f_1,\ldots,f_n)\simalg \mathrm{Expr}_2(f_1,\ldots,f_n),\]
but also about the compatibility of the  operations with algebraic homotopy, i.~e.
\[\forall i\, f_i\simalg  f_i'\quad\Rightarrow\quad
\mathrm{Op}(f_1,\ldots,f_n)\simalg \mathrm{Op}(f_1',\ldots,f_n').\]
Homotopy compatibility is rather trivial in the topological setting but it is less
trivial in the algebraic setting.
Nevertheless, we can reduce this problem:
\begin{itemize}
\item[(i)] We can check compatibility in the variables separately.
\item[(ii)] Even there, it is sufficient to check it in two special cases:
First we must check invariance for stabilization, i.~e.  the case
$f'_i=f\oplus\mathbf 0_{\Xi_0}\oplus\mathbf 0_{\Xi_1}'$,
and, second,  to check conjugation invariance when $f_i\xrightarrow{g}f'_i$.
\end{itemize}
Taking direct sum of operators, we can see easily that
the sum operation is compatible with algebraic homotopy.
We do not have to worry about the homotopy compatibility of the operations $\mathbf 0$ and $\mathbf 1$.
The homotopy compatibility of the additive inverse and the product
are left to be demonstrated, although these problems are not equally hard.
We summarize what are the identities we want to prove:
\begin{enumerate}
\item Addivite inverse.
\item Homotopy compatibility of the inverse.
\item Equality of products (or multiplicative commutativity).
\item Homotopy compatibility of the products.
\end{enumerate}
\end{point}

\section{Tools: regularization and taming}
In this section, we introduce some tools to deal with algebraic homotopy effectively.
Regarding the definition, one might believe that we allow  excessively large classes
of objects and  morphisms in our Grassmannians.
We will show that this is not the case.
By ``regularization'', one can reduce the variety of base objects,
and by ``taming'', one can replace the morphisms by smooth ones.
\begin{lemma}[\textbf{Virtual cancellation}] \label{lem:thelarge}
For $a\in\Psi_{\Omega}(\mathfrak A^+)$, we have $\langle a,a\rangle \simalg \mathbf 0.$
\end{lemma}
\noindent
\textit{Remark.} The exact meaning of this statement is the existence of an algebraic homotopy
of functions $\langle \id_{\invol(\Psi_{\Omega}^{(2)}(\mathfrak A^+))},
\id_{\invol(\Psi_{\Omega}^{(2)}(\mathfrak A^+))}\rangle
\simalg \mathbf 0$ with domain $\invol(\Psi_{\Omega}^{(2)}(\mathfrak A^+))$, but we allow
the colloquiality of using variables instead of functions here, and in the future.
\begin{proof}
In terms of    $(\mathbb Z+\frac12)\pmb\times\mathbb Z$
and $\mathbb Z\pmb\times(\mathbb Z+\frac12)$  block matrices, let us consider
\[\mathsf B(a)=\sum_{n\in\mathbb Z} a \mathbf e_{n-\frac12,n}+\bar a\mathbf e_{n+\frac12,n},
\quad\text{and}\quad
\mathsf B(a)^{-1}=\sum_{n\in\mathbb Z} a \mathbf e_{n,n-\frac12}+\bar a\mathbf e_{n,n+\frac12}. \]
If the index set $\{0\}\times\Omega$ is replaced with $\Omega$, then it yields
\[1_{\mathbb Z^-\times\Omega}\oplus a \oplus 0_{\mathbb Z^+\times\Omega}\xrightarrow{\mathsf B(a)}
1_{(-\frac12-\mathbb N)\times\Omega}\oplus 0_{(\frac12+\mathbb N)\times\Omega}.\]
After doubling the terms, it provides an algebraic homotopy as required.
\end{proof}
\begin{point} \textbf{Regularization.}
For $\langle b,a\rangle\in \mathcal G^{(2)}_{\Omega}(\mathfrak A)$, we define its regularized
matrix as the $\{0,1\}\times\Omega$ block matrix
\[\mathrm R\langle b,a\rangle:=\sw_{01}(a)(b\oplus\bar a)\sw_{01}(a)=
\begin{bmatrix}\bar a(b-a)\bar a&\bar a(b-a)a\\a(b-a)\bar a&1_\Omega+ a(b-a) a\end{bmatrix}.\]
Then
$\mathrm R\langle b,a\rangle\approx\mathrm R\langle a,a\rangle
=0_\Omega\oplus 1_\Omega=\begin{bmatrix}0_\Omega&\\&1_\Omega\end{bmatrix}$,
so we can define the regularized pair
\[\mathbf R\langle b,a\rangle:=\left\langle\mathrm R\langle b,a\rangle,0_\Omega\oplus 1_\Omega\right\rangle.\]
The trivial pair is
\[\mathbf R\mathbf 0_\Omega=\langle0_\Omega\oplus 1_\Omega,0_\Omega\oplus 1_\Omega\rangle.\]
\end{point}
\begin{point} \textbf{Regularization of morphisms.}
Let $\langle\psi,\phi\rangle\in\Psi^{(2)}_{\Omega',\Omega}(\mathfrak A^+)^\star$, and let
$a\in \Psi_\Omega(\mathfrak A^+)$ be an idempotent.
The  regularized morphism is the $\{0,1\}\times\Omega'\pmb \times \{0,1\}\times\Omega$ block matrix
\[\mathrm R\langle\psi,\phi\rangle,a):=\sw_{01}(a{}^\phi)(\psi\oplus\phi)\sw_{01}(a).\]
Then
$\mathrm R(\langle\psi,\phi\rangle,a)\approx \mathrm R(\langle\phi,\phi\rangle,a)=\phi\oplus\phi,$
so we can define
\[\mathbf R(\langle\psi,\phi\rangle,a):=\langle\mathrm R(\langle\psi,\phi\rangle,a),\phi\oplus\phi\rangle.\]
It yields
\[\mathbf R\langle b,a\rangle\xrightarrow{\mathbf R(\langle\psi,\phi\rangle,a)}
\mathbf R\langle b{}^\psi,a{}^\phi\rangle.\]
Based upon this, it also reasonable to write
$\mathbf R(\langle\psi,\phi\rangle_{\langle b,a\rangle})$ instead of
$\mathbf R(\langle\psi,\phi\rangle,a)$, etc. even if there is no dependence on $b$.
\end{point}
\begin{lemma}  $\mathbf R$  is an operation algebraically homotopic to the identity:
\[\langle b,a\rangle\simalg \mathbf R\langle b,a\rangle.\]
\begin{proof} By cancellation and conjugation
$\langle b,a\rangle\simalg \langle b,a\rangle\oplus\langle\bar a,\bar a\rangle
\xrightarrow{\langle\sw_{12}(a),\sw_{12}(a)\rangle} \mathbf R\langle b,a\rangle$.
\end{proof}
\end{lemma}
The transitivity of algebraic homotopy automatically ensures that $\mathbf R$ is compatible with algebraic homotopy.
Nevertheless, one can also show this directly, using the straightforward compatibility with direct sums
and regularized morphisms.
Ultimately,  by regularization, we can bring the base terms into simple form.
Another natural expectation is that morphisms should not deviate much from identity.
This can be achieved as follows.
\begin{point} \textbf{Translation.}
Using $\{0,0',0''\}\times\Omega$ block matrices we set
\[\mathrm H\langle b,a\rangle:=\sw_{0'0''}(a)\sw_{0'0''}(b)\sw_{00'}(b)\sw_{00'}(a).\]
Then
$\mathrm H\langle b,a\rangle\approx\mathrm H\langle a,a\rangle= 1_{\{0,0',0''\}\times\Omega}$,
so we can define
\[\mathbf H\langle b,a\rangle:=\langle\mathrm H\langle b,a\rangle,
1_{\{0,0',0''\}\times\Omega}\rangle.\]
It yields
\[\langle b\oplus\bar a\oplus a,a\oplus\bar a\oplus a\rangle\xleftarrow{\mathbf H\langle b,a\rangle}
\langle a\oplus\bar a\oplus b,a\oplus\bar a\oplus a\rangle.\]
\end{point}
\begin{point} \textbf{Regularized translation.}
The regularized version is given by $\{0,1,0',1',0'',1''\}\times\Omega$ block matrices as follows: Let
\[\mathrm{HR}\langle b,a\rangle:=
\sw_{01}(a)\sw_{0'1'}(\bar a)\sw_{0''1''}(a)
\sw_{0'0''}(a)\sw_{0'0''}(b)\sw_{00'}(b)\sw_{00'}(a)
\sw_{01}(a)\sw_{0'1'}(\bar a)\sw_{0''1''}(a).\]
Then
$\mathrm{HR}\langle b,a\rangle\approx\mathrm{HR}\langle a,a\rangle=1_{\{0,1,0',1',0'',1''\}\times\Omega},$
which allows us to define
\[\mathbf{HR}\langle b,a\rangle:=\langle\mathrm{HR}\langle b,a\rangle,
1_{\{0,1,0',1',0'',1''\}\times\Omega}\rangle .\]
This yields
\[\mathbf R\langle b,a\rangle\oplus\mathbf{R0}_\Omega \oplus\mathbf{R0}_\Omega
\xleftarrow{\mathbf{HR}\langle b,a\rangle}
\mathbf{R0}_\Omega\oplus\mathbf{R0}_\Omega \oplus\mathbf R\langle b,a\rangle.\]
\end{point}
\begin{point} \textbf{Taming.}
If $\langle b,a\rangle\in\mathcal G^{(2)}_\Omega(\mathfrak A)$,
$\langle\psi,\phi\rangle\in \Psi^{(2)}_\Omega(\mathfrak A^+)^\star $, then let
\[\mathrm T\bigl({\langle\psi,\phi\rangle}_{\langle b,a\rangle}\bigr):=
(\psi\oplus1_\Omega\oplus1_\Omega)
\mathrm H\langle b,a\rangle(\phi^{-1}\oplus1_\Omega\oplus1_\Omega)(\mathrm H\langle b,a\rangle)^{-1}\]
One can see that
$\mathrm T\bigl({\langle\psi,\phi\rangle}_{\langle b,a\rangle}\bigr) \approx
\mathrm T\bigl({\langle\phi,\phi\rangle}_{\langle a,a\rangle}\bigr)=1_{\{0,0',0''\}\times\Omega},$
so we can define
\[\mathbf T\bigl({\langle\psi,\phi\rangle}_{\langle b,a\rangle}\bigr):=\langle
\mathrm T\bigl({\langle\psi,\phi\rangle}_{\langle b,a\rangle}\bigr),1_{\{0,0',0''\}\times\Omega}\rangle.\]
If $\phi$ commutes with $a$, then
\[\langle b\oplus\bar a\oplus a,a\oplus\bar a\oplus a\rangle
\xrightarrow{\mathbf T\left({\langle\psi,\phi\rangle}_{\langle b,a\rangle}\right)}
\langle b{}^\psi\oplus\bar a\oplus a,a\oplus\bar a\oplus a\rangle.\]
In particular, using the variant
\[\mathbf T'\bigl({\langle\psi,\phi\rangle}_{\langle b,a\rangle}\bigr)=
\langle\sw_{0'0''}(a)\mathrm T(\langle\psi,\phi\rangle_{\langle b,a\rangle})\sw_{0'0''}(a)
,1_{\{0,0',0''\}\times\Omega}\rangle\]
we obtain the following
\begin{cor}[\textbf{Stable taming}]\label{cor:tame}
$\langle b,a\rangle\xrightarrow{\langle\psi,\phi\rangle}\langle\tilde b,a\rangle$ implies
\[\left\langle b\oplus 1_\Omega\oplus0_\Omega, a\oplus 1_\Omega\oplus0_\Omega \right\rangle
\xrightarrow{\mathbf T'\bigl({\langle\psi,\phi\rangle}_{\langle b,a\rangle}\bigr)}
\langle\tilde b\oplus 1_\Omega\oplus0_\Omega, a\oplus 1_\Omega\oplus0_\Omega\rangle,\]
where the conjugating base term is $1_{\{0,0',0''\}\times\Omega}$.
\end{cor}
This amounts to the statement that, up to stabilization, all algebraic homotopies
can be realized by smooth morphisms whenever they have a chance.
\end{point}
\begin{point} \textbf{Regularized taming.}
The commutation assumption is satisfied automatically if we apply taming after regularization.
In terms of $\{0,1,0',1',0'',1''\}\times\Omega$ block matrices, let
\begin{multline}
\mathrm{TR}\bigl({\langle\psi,\phi\rangle}_{\langle b,a\rangle}\bigr)
:=(\mathrm R(\langle\psi,\phi\rangle,a)\oplus1_{\{0,1\}\times\Omega}\oplus1_{\{0,1\}\times\Omega})
\mathrm{HR}\langle b,a\rangle\cdot\\\cdot(\mathrm R(\langle\phi,\phi\rangle,a)
\oplus1_{\{0,1\}\times\Omega}\oplus1_{\{0,1\}\times\Omega})(\mathrm{HR}\langle b,a\rangle)^{-1}.
\notag\end{multline}
Then
$\mathrm{TR}\bigl({\langle\psi,\phi\rangle}_{\langle b,a\rangle}\bigr)
\approx \mathrm{TR}\bigl({\langle\phi,\phi\rangle}_{\langle a,a\rangle}\bigr)
=1_{\{0,1,0',1',0'',1''\}\times\Omega}$, so we can set
\[\mathbf{TR}\bigl({\langle\psi,\phi\rangle}_{\langle b,a\rangle}\bigr):=
\langle\mathrm{TR}\bigl({\langle\psi,\phi\rangle}_{\langle b,a\rangle}\bigr)
,1_{\{0,1,0',1',0'',1''\}\times\Omega}\rangle.\]
It yields
\[\mathbf R\langle b,a\rangle\oplus\mathbf{R0}_\Omega\oplus\mathbf{R0}_\Omega
\xrightarrow{  \mathbf{TR}\left({\langle\psi,\phi\rangle}_{\langle b,a\rangle}\right)}
\mathbf R\langle b{}^\psi,a{}^\phi\rangle\oplus\mathbf{R0}_\Omega\oplus\mathbf{R0}_\Omega.\]
\end{point}
\begin{point}
\textbf{Regular relabeling.}
Let us mention one last operation homotopic to the identity, which, however, depends on the index set.
Suppose that $\theta:\Omega\rightarrow \mathbb N$ is a relabeling map.
Then consider the operation $\mathbf R\theta_*$, which is the composition of relabeling
in both components followed by regularization.
From the preceding discussion, it is clear that the map
\[\mathbf R\theta_*:\mathcal G^{(2)}_\Omega(\mathfrak A)\rightarrow \mathcal
G^{(2)}_{\{0,1\}\times\mathbb N}(\mathfrak A)\]
\[\langle b,a\rangle\mapsto\left\langle \mathrm R\langle\theta_*b,\theta_*a \rangle,\mathrm R
\mathbf0_{\mathbb N}
\right\rangle\]
is homotopic to the identity. In fact, the target can be considered to be the smaller space
\[\mathcal G_{\mathbb N}(\mathfrak A):=\{\langle b,\mathrm R\mathbf 0_{\mathbb N}
\rangle\,:b\in\Psi_{\{0,1\}\times\mathbb N}(\mathfrak A^+),\,b^2=b, b\approx \mathrm R\mathbf 0_{\mathbb N} \},\]
the standard infinite single-space Grassmannian.
\end{point}
\begin{point}\label{po:polyn}
\textbf{On polynomial constructions}
We may observe that our constructions are all finite matrix polynomials in terms of
the initial data, except the conjugating matrices $\mathsf B(a)$ in the proof of Lemma \ref{lem:thelarge}.
But even that term was a matrix of finite Toeplitz type. However, even conjugating
matrices of finitely multiple Toeplitz type will become essentially finite matrices if we apply the
taming construction, because the finitely many components of diagonal type will cancel out,
leaving an infinite  but inert part of identity on stabilization terms.
As a consequence, if $\langle b,a \rangle,\langle\tilde b,\tilde a\rangle
\in\mathcal G^{(2)}_{\Omega}(\mathfrak A)$ and
\[\mathbf R\langle b,a \rangle\simalg  \mathbf R\langle\tilde b,\tilde a\rangle\]
such that $\simalg $ is realized by a  of finitely multiple Toeplitz type block matrix, then
due to the universally applicable regularized taming construction
\[\mathbf R\langle b,a \rangle
\simalg
\mathbf R\mathbf R\langle b,a \rangle
\simalg
\mathbf R\mathbf R\langle\tilde b,\tilde a\rangle
\simalg
\mathbf R\langle\tilde b,\tilde a\rangle,\]
where  all $\simalg $ are realized by  finite block matrices.
That amounts to the fact, which will extend to our later experience,
and which can also be checked case by case, that the regularized operations and their
related algebraic homotopies
can be realized by finite block matrices in terms of the initial data.
\end{point}

\section{Establishing the algebraic $H$-ring structure}
\begin{point}
\textbf{Rule of the additive inverse.} In terms of $\{0,0'\}\times\Omega$ block matrices, we have
\[\langle b,a\rangle_\Omega\oplus \langle b,a\rangle_\Omega^{\mathrm{inv}}
\xrightarrow{\langle\sw_{00'}(b),\sw_{00'}(a)\rangle}
\langle0_\Omega\oplus 1_\Omega,0_\Omega\oplus 1_\Omega\rangle\simalg \mathbf 0.\]
\end{point}
\begin{point}
\textbf{Compatibility of the additive inverse.}
The stabilization part follows from
\[(\langle b,a\rangle\oplus\mathbf 0_{\Omega_0}\oplus \mathbf 0'_{\Omega_1})^{\mathrm{inv}}=
\langle\bar b,\bar a\rangle\oplus\mathbf 0'_{\Omega_0}\oplus \mathbf 0_{\Omega_1}
\simalg \langle\bar b,\bar a\rangle=\langle b,a\rangle^{\mathrm{inv}}.\]
The conjugation part is obvious from
$\langle b,a\rangle^{\mathrm{inv}}\xrightarrow{\langle\psi,\phi\rangle}
\langle b{}^\psi,a{}^\phi\rangle^{\mathrm{inv}}.$
\end{point}
\begin{point}\textbf{Commutativity of the product.}
Using $\{0,1,0',1'\}\times\Omega\times\Xi$ matrices,  we set
\[\stackrel{\hookrightarrow}{\mathrm C}(\langle b,a\rangle,\langle d,c\rangle):=
\sw_{10}(b \otimes c+a \otimes \bar c)\sw_{10'}(\bar a \otimes \bar c)
\sw_{10}(\overline{b \otimes d})\sw_{1'0'}(\bar a \otimes \bar c)
\sw_{10'}(\bar a \otimes \bar c)\sw_{10}(a \otimes  d+\bar a \otimes c)\]
One can see that
\[\stackrel{\hookrightarrow}{\mathrm C}(\langle b,a\rangle,\langle d,c\rangle)\approx
\stackrel{\hookrightarrow}{\mathrm C}(\langle a,a\rangle,\langle c,c\rangle) =
\sw_{00'}(\bar a\otimes \bar c)\sw_{11'}(\bar a \otimes\bar c), \]
hence one can define
\[\stackrel{\hookrightarrow}{\mathbf C}(\langle b,a\rangle,\langle d,c\rangle):=
\langle\stackrel{\hookrightarrow}{\mathrm C}(\langle b,a\rangle,\langle d,c\rangle),
\sw_{00'}(\bar a\otimes\bar c)\sw_{11'}(\bar a\otimes\bar c)\rangle.\]
Then, computation yields that
\[\mathbf R(\langle b,a\rangle\overleftarrow\otimes\langle d,c\rangle)\oplus\mathbf{R0}_{\Omega\times\Xi}
\xrightarrow{\stackrel{\hookrightarrow}{\mathbf C}(\langle b,a\rangle,\langle d,c\rangle),}
\mathbf R(\langle b,a\rangle\overrightarrow\otimes\langle d,c\rangle)\oplus\mathbf{R0}_{\Omega\times\Xi},\]
which demonstrates commutativity.
\end{point}
\begin{point} \textbf{Compatibility of the product.}
First, we  prove the compatibility of $\overleftarrow\otimes$ in the second variable.
The stabilization part follows from distributivity and
\[\langle b,a\rangle\overleftarrow\otimes\langle1_{\Xi_1}\oplus0_{\Xi_0},1_{\Xi_1}\oplus0_{\Xi_0}\rangle=
\langle1_{\Omega\times\Xi_1}\oplus0_{\Omega\times\Xi_0},1_{\Omega\times\Xi_1}\oplus0_{\Omega\times\Xi_0}\rangle.\]
The conjugation part follows from
\[\langle b,a\rangle\overleftarrow\otimes\langle d,c\rangle
\xrightarrow{\langle b\otimes\theta+\bar b\otimes \chi,a\otimes\theta+\bar a\otimes\chi\rangle}
\langle b,a\rangle\overleftarrow\otimes\langle d{}^\theta,c{}^\chi\rangle.\]
One can do the compatibility of $\overrightarrow\otimes$ in the first variable in the same way.
But then, the equivalence of $\overleftarrow\otimes$ and $\overrightarrow\otimes$,
from the previous point, implies compatibility in each variable.
\end{point}
\begin{remark}
 One might prefer the variant additive inverse, the classical switch operation
\[\langle b,a\rangle^{\mathrm{inv}\prime}=\langle a,b\rangle.\]
In order to show equivalence to the usual additive inverse, it is sufficient to show that the variant operation
\[\langle b,a\rangle'=\langle\bar a,\bar b\rangle\]
is algebraically homotopic to the identity operation. But this follows from
\[\mathbf R\langle a,b\rangle
\xrightarrow{\langle\sw_{01}(b)\sw_{01}(a),1_{\{0,1\}\times\Omega}\rangle}
\mathbf R\langle a,b\rangle'.\]
\end{remark}
\begin{cor}\label{cor:diff} For idempotents $b\approx\tilde b\approx a$, we have
$\langle b,\tilde b\rangle\simalg
\langle b,a\rangle\oplus \langle\tilde b,a \rangle^{\inv}.$
\begin{proof} Indeed, we see
$\langle b,\tilde b\rangle\simalg \langle a\oplus\bar a\oplus b,a\oplus\bar a\oplus\tilde b\rangle
\xrightarrow{\mathbf H\langle b,a\rangle}\langle b\oplus\bar a\oplus a,a\oplus\bar a\oplus\tilde b\rangle
\simalg \langle b,a\rangle\oplus\langle a,\tilde b\rangle\simalg
\langle b,a\rangle\oplus\langle\tilde b,a\rangle^{\inv}$.
\end{proof}
\end{cor}
\begin{remark}
In  classical algebraic $K$-theory,
the class of $\langle b,a\rangle\overleftarrow\otimes \langle d,c\rangle$ is
$[b\otimes d +\bar b\otimes c]-[a\otimes d +\bar a\otimes c]=
[b][d]+([1_\Omega]-[b])[c]-[a][d]-([1_\Omega]-[a])[c]=([b]-[a])([d]-[c])$, i.~e. the product of the
classes of $\langle b,a\rangle $ and $\langle d,c\rangle$.
\end{remark}
\begin{point} \textbf{Summary.}
We have proved that the functor
\[\Omega,\mathfrak A\mapsto \mathcal G_\Omega^{(2)}(\mathfrak A)\]
is a commutative unital algebraic homotopy ring.
This is a very flexible construction, but it has some strange features:
\begin{itemize}
\item Addition changes the index sets $\Omega_1,\Omega_2\mapsto \Omega_1\dot\cup\Omega_2.$
\item Multiplication also changes the index sets
$\Omega_1,\Omega_2\mapsto \Omega_1\times\Omega_2$, but it also changes the rings
$\mathfrak A,\mathfrak B\rightarrow\mathfrak A\otimes\mathfrak B$. We have certain flexibility in
what kind of tensor products we consider.
\item The multiplicative unit element calls for a base ring like $\mathbb Z$ or $\mathbb R$,
which may not be embedded into $\mathfrak A$.
\end{itemize}
Our ultimate objective is, however, to study single-space Grassmannians.
\end{point}
\begin{point}\textbf{Single-space Grassmannian.} Suppose that $\mathfrak A$ is a commutative ring.
(The multiplication can be considered as tensor product.)
Consider the space $\mathcal G_{\mathbb N}(\mathfrak A)$. Then we can consider the alternative oprations
\[a\tilde+b:=\mathbf R\theta_1{}_*(a\oplus b),\quad
\tilde -a:=\mathbf R\theta_2{}_*(a^{\mathrm{inv}}),\quad
\tilde 0:=\mathbf{R0}_{\mathbb N},\quad
a\tilde\otimes b=\mathbf R\theta_3{}_*(a\otimes b),\]
and, if $\mathfrak A$ is unital
\[ \tilde1=\mathbf R\theta_4{}_*\mathbf 1;\]
where the maps
$\theta_1:\{0,1\}\times \mathbb N \dot\cup  \{0,1\}\times \mathbb N\rightarrow \mathbb N,$
$\theta_2:\{0,1\}\times \mathbb N\rightarrow \mathbb N,$
$\theta_3:\{0,1\}\times \mathbb N\times \{0,1\}\times\mathbb N\rightarrow \mathbb N,$
$\theta_4:\{0\}\rightarrow \mathbb N$
are some fixed, not necessarily isomorphic relabeling maps.
These new operations are algebraic operations in strict sense.
From the algebraic homotopy ring structure of the functorial Grassmannian, it
is immediately clear that the new operations yield an algebraic homotopy ring by restriction.

Moreover, due to the taming construction described earlier, we can use
algebraic homotopy in strong sense: we say that
$f_1:X\rightarrow \mathcal G_{\mathbb N}(\mathfrak A)$
and $f_2:X\rightarrow \mathcal G_{\mathbb N}(\mathfrak A)$ are smoothly algebraically
homotopic
\[f_1\simeq_{\mathrm{sm\, alg}}f_2\]
if there is a map $g:X\rightarrow \mathcal K_{\{0,1\}\times\{0,1\}\times\mathbb N}(\mathfrak A)^\star$ such that
\[f_2\oplus\mathbf R0_{\mathbb N}\xrightarrow g f_1 \oplus\mathbf R0_{\mathbb N} .\]
It is an equivalence relation. Transitivity, for example, follows as in the general case, except
one needs a final relabeling in the  stabilizing indices.

Furthermore, as the operations are sufficiently regularized, and as it was explained \ref{po:polyn},
the operations and the realizing algebraic homotopies  can be chosen to be matrix
polynomials. Here  we take the extended sense that tensor products and
relabeling matrices which are constant are still in the notion of polynomiality.
Hence, we obtain:
\end{point}
\begin{theorem}
If $\mathfrak A$ is a (unital) commutative polymetric ring, then
$\mathcal G_{\mathbb N}(\mathfrak A)$ with the operations above
is a (unital) commutative ``strong polynomial'' algebraic $H$-ring.
By ``strong polynomial'' we mean that the algebraic homotopies
(for the identities and compatibilities)
are induced by maps
$\mathcal G_{\mathbb N}(\mathfrak A)\times\ldots\times\mathcal G_{\mathbb N}(\mathfrak A)
\rightarrow\mathcal K_{\{0,1\}\times\{0,1\}\times\mathbb N}(\mathfrak A)^\star$
assembled polynomially.
\end{theorem}
\begin{point}
It must be clear that the operations on the Grassmannian allow
several variants, but they are OK as long as one can deduce the existence of
algebraic homotopies from the functorial Grassmannian.
It must be also clear that  several intermediate
constructions are allowed, like a product
$\mathcal G_{\mathbb N}(\mathfrak A)\times\mathcal G_{\mathbb N}(\mathfrak B)
\rightarrow\mathcal G_{\mathbb N}(\mathfrak A\otimes\mathfrak B)$.
This comment also applies for later constructions, but we will not emphasize it further.
\end{point}
\section{Establishing the smooth topological $H$-ring structure}
This section applies to locally convex algebras $\mathfrak A$.
The key point is that we can deal with stabilization internally, without adjoining extra variables.
In short terms, we can make extra space using homotopies.
\begin{point}
First, we discuss the stabilization
of the  algebra $\mathcal K_{\mathbb N}(\mathfrak A)$, which will be extended to the larger idempotents.
Stabilization can be organized as in \cite{SS}.
For $\theta\in[0,\frac\pi2]$  we consider
\[\mathsf C(\theta)=\begin{bmatrix}
s&ts&t^2s&t^3s&t^4s&\cdots\\
-t & s^2& ts^2 & t^2s^2  & t^3s^2&\ddots\\
&-t & s^2 & ts^2 & t^2s^2 &\ddots\\
&&-t & s^2 & ts^2 &\ddots\\
&&&-t & s^2 &\ddots\\
&&&&-t &\ddots\\
&&&&&\ddots\\
\end{bmatrix},\]
where $t=\sin\theta$ and $ s=\cos\theta$.
Then a stabilizing homotopy is given by
\[T_{\mathcal K}:\mathcal K_{\mathbb N}(\mathfrak A)\times \left[0,\frac\pi2\right]
\rightarrow\mathcal K_{\mathbb N}(\mathfrak A)\]
\[A,\theta\mapsto T_{\mathcal K}(A,\theta)=\mathsf C(\theta)A\mathsf C(\theta)^\top.\]

It yields a smooth homotopy between the identity
$T_{\mathcal K}{}(\cdot,0)=\id_{\mathcal K_{\mathbb N}(\mathfrak A)}$
and the relabeling map
$T_{\mathcal K}(\cdot,\frac\pi2)=r_*$,
where $r(n)=n+1$. It stabilizes by one extra matrix entry, which may seem insufficient.
However, using it as an $\mathbb N\times\Omega$ block matrix construction, it achieves
stabilization by infinitely many entries.
\end{point}
\begin{point}
There is a slightly more complicated, yet elegant way to achieve stabilization:
We can consider $\mathsf C(\theta)$ as stabilization of the space
$\mathcal H_{\mathbb N}(\mathfrak A)$, and then we consider the odd or even ``quantized''
situation. Without going into details, the topology of $\mathcal H_{\mathbb N}(\mathfrak A)$ is induced by the
seminorms
\[p_\alpha\left(\sum_{n\in\mathbb N}{a_n\mathbf e_n}\right)=\sum_{n\in\mathbb N} p(a_n)\alpha^n
\qquad\qquad(\alpha>0).\]

The odd quantized space $\mathrm{Qu}^1\mathcal H_{\mathbb N}(\mathfrak A) $ is the space with basis
$\mathbf e_{i_1}\wedge \mathbf e_{i_2}\wedge\ldots \wedge\mathbf e_{i_k}$,
where $i_1<i_2<\ldots<i_k$. The seminorms induced ultimately are
\[p\left(\sum a_{i_1,\ldots,i_k}\mathbf e_{i_1}\wedge
\mathbf e_{i_2}\wedge\ldots \wedge\mathbf e_{i_k}\right)=
\sum (\alpha^{i_1}+\ldots+\alpha^{i_k}  )p(a_{i_1,\ldots,i_k}).\]
This quantized  space is isomorphic to the space of rapidly decreasing sequences
$\mathcal S_{\mathbb N}(\mathfrak A)$ by
\[\mathrm V: \qquad\mathbf e_{i_1}\wedge \mathbf e_{i_2}\wedge\ldots \wedge\mathbf e_{i_k}\mapsto
\mathbf e_{2^{i_1}+\ldots+2^{i_k}}.\]

The even quantized space $\mathrm{Qu}^0\mathcal H_{\mathbb N}(\mathfrak A)$
is the space with basis $\mathbf e_{i_1}\odot \mathbf e_{i_2}\odot\ldots \odot\mathbf e_{i_k}$, where
$i_1\leq i_2\leq\ldots\leq i_k$. The induced seminorms are, similarly,
\[p\left(\sum a_{i_1,\ldots,i_k}\mathbf e_{i_1}\odot
\mathbf e_{i_2}\odot\ldots \odot\mathbf e_{i_k}\right)=
\sum (\alpha^{i_1}+\ldots+\alpha^{i_k}  )p(a_{i_1,\ldots,i_k}).\]
This is also isomorphic to $\mathcal S_{\mathbb N}(\mathfrak A)$, although this is
less transparent. For this reason, it is more practical to use the odd quantization.

We can consider the quantized matrices
$\mathrm{Qu}^0\mathsf C(\theta)$ and $\mathrm{Qu}^1\mathsf C(\theta)$.
For example,
\[\mathrm{Qu}^1\mathsf C(\theta)(\mathbf e_0\wedge\mathbf e_1)=
(s\mathbf e_0-t\mathbf e_1)\wedge(st\mathbf e_0+s^2\mathbf e_1-t\mathbf e_2)=
s\mathbf e_0\wedge\mathbf e_1-st\mathbf e_0\wedge\mathbf e_2+t^2\mathbf e_1\wedge\mathbf e_2,\]
or, in other terms,
\[\mathrm V\mathrm{Qu}^1\mathsf C(\theta)\mathbf e_3=s\mathbf e_3-st\mathbf e_5+t^2\mathbf e_6.\]

This yields a smooth homomorphic homotopy
$\mathrm V\mathrm{Qu}^1T_{\mathcal K}$
between
$\mathrm V\mathrm{Qu}^1T_{\mathcal K}(\cdot,0)=\id_{\mathcal K_{\mathbb N}(\mathfrak A)}$
and the ``halving'' relabeling map
$\mathrm V\mathrm{Qu}^1T_{\mathcal K}(\cdot,\tfrac\pi2)=\mathrm{hv}_*$
where $\mathrm{hv}(n)=2n$.
It has the nice property that it achieves infinite stabilization immediately,
and it allows  plenty of direct sum decompositions. In fact, we can consider any homomorphic
homotopy $\mathrm{Hv}$ instead of $\mathrm V\mathrm{Qu}^1T_{\mathcal K}$ as long as it
yields a smooth homotopy as above.
\end{point}
\begin{point}
It is natural to make an identification
\[\mathcal K_{\mathbb N}(\mathfrak A)\equiv
\begin{bmatrix} \mathcal K_{\mathbb N}(\mathfrak A)&\mathcal K_{\mathbb N}(\mathfrak A)
\\\mathcal K_{\mathbb N}(\mathfrak A)&\mathcal K_{\mathbb N}(\mathfrak A)\end{bmatrix}^{\sim}.\]
In what follows, ${}^\sim$   indicates  that we
consider a decomposition of the index set $\mathbb N$ into two copies of $\mathbb N$ through
 the indexing convention
$(0,n)\leftrightarrow 2n,(1,n)\leftrightarrow 1+2n$.
This looks especially simple if we use binary numbers with digits written in reverse order.

$\mathrm{Hv}$ yields a homotopy between
\[\id_{\mathcal K_{\mathbb N}(\mathfrak A)}\qquad\text{and}\qquad
 \id_{\mathcal K_{\mathbb N}(\mathfrak A)}\tilde\oplus 0_{\mathbb N}\]
through algebra homomorphisms.
But we want to stabilize idempotents. This can be done by setting
\[\mathrm{Hv}\left(
\begin{bmatrix}b_{00}&b_{01}\\b_{10}& 1_{\mathbb N}+b_{11} \end{bmatrix}\right):=
\begin{bmatrix}\mathrm{Hv}(b_{00})&\mathrm{Hv}(b_{01})\\
\mathrm{Hv}(b_{10})&1_{\mathbb N}+\mathrm{Hv}(b_{11}) \end{bmatrix}.\]
This yields a smooth homotopy between
\[b\equiv\begin{bmatrix}b_{00}&b_{01}\\b_{10}& 1_{\mathbb N}+b_{11} \end{bmatrix}
\qquad\text{and}\qquad
\begin{bmatrix}b_{00}&b_{01}\\b_{10}& 1_{\mathbb N}+b_{11} \end{bmatrix}\tilde \oplus
\begin{bmatrix}0_{\mathbb N}&\\& 1_{\mathbb N} \end{bmatrix}\equiv
\begin{bmatrix}b&\\& \mathrm R\mathbf 0_{\mathbb N} \end{bmatrix}^{\sim}.\]
In simple terms, this means that by using a smooth homotopy, we can
always make free space, i.~e. we can stabilize.
It also means that we can conjugate freely by smooth terms. Indeed,  using the conjugating matrices
\[\begin{bmatrix}\cos\alpha
&-\hat{\mathrm{hv}}^\top  \sin\alpha\\\hat{\mathrm{hv}}\sin\alpha&1_{2\mathbb N }\cos\alpha+1_{2\mathbb N+1 }
\end{bmatrix}
\begin{bmatrix}\phi&\\&1\end{bmatrix}
\begin{bmatrix}\cos\alpha
&\hat{\mathrm{hv}}^\top  \sin\alpha\\-\hat{\mathrm{hv}}\sin\alpha&1_{2\mathbb N }\cos\alpha+1_{2\mathbb N+1 }
\end{bmatrix},\]
we quickly obtain a smooth homotopy between
\[\begin{bmatrix}b&\\& \mathrm R\mathbf 0_{\mathbb N} \end{bmatrix}^{\sim}
\quad\text{and}\quad
\begin{bmatrix}b{}^\phi&\\& \mathrm R\mathbf 0_{\mathbb N} \end{bmatrix}^{\sim}.\]

As the homotopies used are smooth, and the conjugating terms used are smooth
(in fact polynomials) in terms of the initial data, we obtain
\begin{theorem}
If $\mathfrak A$ is a (unital) commutative locally convex algebra, then
$\mathcal G_{\mathbb N}(\mathfrak A)$
is a (unital) commutative ``smooth'' topological $H$-ring.
By ``smooth'' we mean that the topological homotopies are induced by smooth maps
$\mathcal G_{\mathbb N}(\mathfrak A)\times\ldots\times\mathcal G_{\mathbb N}(\mathfrak A)\times[0,1]
\rightarrow\mathcal G_{\mathbb N}(\mathfrak A).$
\end{theorem}

Now, the exact meaning of ``smooth'' in the setting of infinite-dimensional
spaces remains somewhat unclear, but if one tries his favorite
definition,  then he will most likely agree.
\end{point}

\section{Establishing the smooth topological algebraic $H$-ring structure}
\begin{point}
We say that the polymetric ring $\mathfrak A$ is  strong if for every seminorm $p$ there is a seminorm
$\tilde p$ such that
$p(a_1\ldots a_n)\leq \tilde p(a_1)\ldots \tilde p(a_n)$
for any $n\in\mathbb N$. These algebras behave well with respect to forming $\mathfrak A^+$,
and  $\mathcal K_\Omega(\mathfrak A)$, etc.

 Suppose that $\mathfrak A$ is a strong locally convex algebra.
Let $[a,b]\subset\mathbb R$ be a closed  interval. If $A:[a,b]\rightarrow\mathfrak A$ is a, say, continous
function, then $C(t)=A(t)\,\mathrm dt$ is a continuous ordered measure on $[a,b]$.
For such a continuous measure, the time-ordered exponential
\[\exp_t C(t)=1+\int_{t_1}C(t_1)+\ldots+ \int_{t_1\leq\ldots\leq t_n}
C(t_n)\ldots C(t_1)+\ldots \]
can be considered. One can check
\end{point}
\begin{lemma} If $\mathfrak A$ is a strong locally convex algebra, and
 $P:[a,b]\rightarrow\mathfrak A$ is a smooth idempotent-valued map, then the map
\[A_p:\{(t_1,t_2)\,:\,t_1,t_2\in[a,b], t_1\leq t_2\}\rightarrow\mathfrak A\]
\[(t_1,t_2)\mapsto A_P(t_1,t_2)=\exp_t \dot P(t)P(t)-P(t)\dot P(t)\,\mathrm dt|_{[t_1,t_2]}\]
is also smooth, and
\[P(t_1)\xrightarrow{A_P(t_1,t_2)} P(t_2).\]
\end{lemma}
\begin{remark}
What happens above corresponds to parallel transport along a connection, which can be written in local form as
$\nabla=\mathrm d- \mathrm dP\,P+P\,\mathrm dP= P.\mathrm d.P+ (1-P).\mathrm d.(1-P)\,.$
\end{remark}
Consequently, smooth homotopies can be lifted to conjugation. Hence we obtain
\begin{theorem}
If $\mathfrak A$ is a (unital) commutative strong locally convex algebra, then
$\mathcal G_{\mathbb N}(\mathfrak A)$
is a (unital) commutative ``smooth topological algebraic'' $H$-ring.
By ``smooth topological algebraic'' we mean that the topological homotopies are induced by conjugating
by smooth maps
$\mathcal G_{\mathbb N}(\mathfrak A)\times\ldots\times\mathcal G_{\mathbb N}(\mathfrak A)\times[0,1]
\rightarrow\mathcal K_{\{0,1\}\times\mathbb N}(\mathfrak A)^\star.$
\end{theorem}
\section{The core information of virtual idempotents}
\begin{point}
Let us use the abbreviations
$h_{ij}:=a^{[i]} (b-a)a^{[j]}$,
where $a^{[0]}=a$ and $ a^{[1]}=\bar a $. Then
\[\mathrm R(\langle b,a\rangle)=\begin{bmatrix}h_{00}&h_{01}\\h_{10}&1_\Omega+h_{11}\end{bmatrix}.\]
In other terms, the $h_{ij}$'s contain the information what is left after regularization.
It turns out that the regularization of various operations can be expressed in terms of that data:
\[\mathrm R(\langle b,a\rangle^{\mathrm{inv}})=
\begin{bmatrix}h_{11}&h_{10}\\h_{01}&1_\Omega+h_{00}\end{bmatrix};\]
and similarly, for
\[\mathrm R(\langle b,a\rangle')=\begin{bmatrix}h'_{00}&h'_{01}\\h'_{10}&1_\Omega+h'_{11}\end{bmatrix},\]
it yields
\[=\sw_{01}( h_{00}-h_{01}+h_{10}-h_{11})
\begin{bmatrix}h_{00}&h_{01}\\h_{10}&1_\Omega+h_{11}\end{bmatrix}
\sw_{01}(h_{00}+h_{01}-h_{10}-h_{11})\]
\[=\begin{bmatrix}h_{00}&h_{10}\\h_{01}&1_\Omega+h_{11}\end{bmatrix}-
\bar{\mathrm s}_{01}(h_{01}h_{10}+h_{01}h_{11}+h_{11}h_{10}-h_{10}h_{01}).\]

Furthermore, if we use the abbreviations
\[\mathrm R(\langle d,c\rangle)=\begin{bmatrix}k_{00}&k_{01}\\k_{10}&1_\Omega+k_{11}\end{bmatrix},\qquad
\mathrm R(\langle d,c\rangle')=\begin{bmatrix}k'_{00}&k'_{01}\\k'_{10}&1_\Omega+k'_{11}\end{bmatrix},\]
then it yields
\[\mathrm R(\langle b,a\rangle\overleftarrow{\otimes} \langle d,c\rangle)=
\begin{bmatrix}0_{\Omega\times\Xi}&\\&1_{\Omega\times\Xi}\end{bmatrix}
+\begin{bmatrix}h_{11}\ootimes k_{00}'&h_{11}\ootimes k_{01}'\\h_{11}\ootimes k_{10}'&h_{11}\ootimes k_{11}'
\end{bmatrix}
+\begin{bmatrix}h_{00}\ootimes k_{00}&h_{00}\ootimes k_{01}\\h_{00}\ootimes k_{10}&h_{00}\ootimes k_{11}
\end{bmatrix}+\]\[
+\begin{bmatrix}&h_{10}\ootimes k_{11}{+}h_{10}\ootimes k_{01}{+}h_{01}\ootimes k_{01}{+}h_{01}\ootimes k_{00}\\
h_{01}\ootimes k_{11}{+}h_{01}\ootimes k_{10}{+}h_{10}\ootimes k_{10}{+}h_{10}\ootimes k_{00}&
\end{bmatrix};\]
and
we obtain a similar formula for $\mathrm R(\langle b,a\rangle\overrightarrow{\otimes} \langle d,c\rangle)$
except the terms $h_{ij}\otimes k_{lm}^{(\prime)}$ should be replaced by
$h_{lm}^{(\prime)}\otimes k_{ij}$.
The additive neutral element corresponds to $\mathrm R\mathbf 0=\mathbf 0$.
The multiplicative neutral idempotent corresponds to the
``scalar'' matrix $\mathrm R\mathbf 1=\begin{bmatrix}1&\\&1\end{bmatrix}$.
\end{point}

\begin{point}
More generally, we can talk about regular idempotents. The general matrix
\[H=\begin{bmatrix}h_{00}&h_{01}\\h_{10}&1_\Omega+h_{11}\end{bmatrix}\]
is an idempotent if and only if
\begin{align}
h_{00}h_{00}+h_{01}h_{10}&=h_{00},&h_{00}h_{01}+h_{01}h_{11}&=0,\notag\\
h_{10}h_{00}+h_{11}h_{10}&=0,&h_{10}h_{01}+h_{11}h_{11}&=-h_{11} \notag.
\end{align}
We say that such an idempotent is  regular  if the identities
\begin{align}
h_{00}h_{10}=h_{01}h_{00}&=0,&h_{00}h_{11}=h_{01}h_{01}&=0,\notag\\
h_{10}h_{10}=h_{11}h_{00}&=0,&h_{10}h_{11}=h_{11}h_{01}&=0.\notag
\end{align}
also hold.
Inspired by the formulas of the previous point, one can define ``regular'' operations
$\oplus,\mathrm{\inv},{}',\overleftarrow\otimes,\overrightarrow\otimes$
for regular idempotents. They satisfy the same identities as their ordinary counterparts;
although checking that they yield new regular idempotents is quite tedious already.
Using the standard machinery, one can show that the ``regular'' operations
are homotopic to the ``ordinary'' operations.
In those terms, the regularization operation $\mathbf R$ is a natural trivial
homomorphism from pairs of idempotents to regular idempotents.
\end{point}
\begin{point}
\textbf{Regular dimension.} Based upon the observations above, we
propose a notion of dimension of virtual idempotents which behaves well with
respect to non-unital rings \textit{and} products:
For $\langle b,a\rangle$ let $\dim\, \langle b,a\rangle$ be the infinum of the cardinality
of those $\Xi$ such that
\[\langle b,a\rangle\simalg \langle U, \mathrm R0_{\Xi}\rangle,\]
where $U$ is a regular idempotent. Then we obtain
\end{point}
\begin{theorem} The dimension defined above has the following properties:

(a) $\mathbf p_1\simalg  \mathbf p_2$ implies $\dim\mathbf  p_1=\dim\mathbf  p_2$.

(b) $\dim \mathbf p=0$ holds if and only if $\mathbf p\simeq \mathbf 0$.

(c) $\dim\mathbf p_1\oplus\mathbf p_2\leq \dim\mathbf p_1+\dim\mathbf p_2 $.

(d) $\dim\mathbf p^{\mathrm{inv}}=\dim\mathbf p.$

(e) $\dim\mathbf p_1 \otimes\mathbf p_2\leq \dim\mathbf p_1\,\dim\mathbf p_2 $.

(f)  $\dim\mathbf 1=1 $, if it applies.

\end{theorem}
\section{Finite dimensionality}
From the previous discussion, it is clear that a Grassmannian element is
finite dimensional if it can be represented by a pair of finite matrices.
If $\mathfrak A$ is a discrete ring, then finite dimensionality holds automatically.
More generally, every Grassmannian element is finite dimensional  if
we can approximate in $\mathfrak A$ in traditional sense.
We demonstrate this statement in order to connect to the classical viewpoint.
\begin{point}
We say that a polymetric $\mathfrak A$ is norm-strong if there
is a distinguished seminorm $q$ such that for each seminorm $p$ there is
an other seminorm $\tilde p$ and $C,i,j\in\mathbb N$ such that
for all $X_1,X_2,\ldots, X_n$ $(n\geq 1)$, we have
\[p(X_1\ldots X_n)\leq Cn^i\sum_{{h:\{1,\ldots, n\}\rightarrow\{q,\tilde p\},\,
\#\{r\,:\,h_r\neq q\}\leq j}} h_1(X_1)\ldots h_n(X_n).\]
This class of polymetric rings also behaves well with respect to taking $\mathfrak A^+$
and $\mathcal K_{\Omega}(\mathfrak A)$.
\end{point}
\begin{point}
Suppose now that $\mathfrak A$ is a norm-strong polymetric ring.
If $P\approx\mathrm R\mathbf 0_\Omega$ is an idempotent, then we can take $\varepsilon$ such that
\[\varepsilon=P-\mathrm R\mathbf 0_\Omega\text{ except in finitely many entries,}  \]
yet so that in a distinguished seminorm $\|\varepsilon\|$ can be arbitrarily small.
Then
\[\widetilde P_\varepsilon=P-\varepsilon\]
differs from $\mathrm R\mathbf 0_\Omega$ only in finitely many entries; however, it is
not necessarily an idempotent.

In general, if $\widetilde P$ is not an idempotent but quite close to being one,
then we may try  the idempotent operation
\[\idem \widetilde P=\int \frac{\widetilde Pz}{(1-\widetilde P)+\widetilde Pz}\frac{|dz|}{2\pi},\]
where we apply formal integration for rapidly decreasing Laurent series.
For the sake of brevity, we will write $1$ for $1_{\{0,1\}\times\Omega}$ in the rest of the section.

In the present case, due to the norm-strong structure,  this machinery works using classical Neumann series.
Indeed, if $\varepsilon$ is so small that
$\|P\varepsilon\|+\|(1-P)\varepsilon\|<\tfrac12$
then
\[P_\varepsilon:=\idem \widetilde P_\varepsilon=\int
(P-\varepsilon)z ((1-P)+Pz^{-1})\Bigl(1-\varepsilon(z-1) ((1-P)+Pz^{-1})\Bigr)^{-1}\frac{|dz|}{2\pi}\]
\[=P{+}\left(2P\varepsilon P{-}\varepsilon P{-}P\varepsilon\right){+}
\left(P\varepsilon\varepsilon{+}\varepsilon P\varepsilon{+}\varepsilon \varepsilon P
{-}3P\varepsilon P\varepsilon {-}3P\varepsilon\varepsilon P{-}3\varepsilon P\varepsilon P
{+}6P\varepsilon P\varepsilon P\right){+}\ldots\]
is convergent.
According to its definition, $\idem \widetilde P_\varepsilon$ is an idempotent which differs from
$\mathrm R\mathbf 0_\Omega$ in finitely many terms.
It is easily to quantify that if  $\varepsilon$ is very small, then $\idem \widetilde P_\varepsilon$
is very close to $P$.

Now, if an idempotent $Q$ is sufficiently close to the idempotent $P$, then
$1-P-Q=(1-2Q)(1-P-Q+2QP)$
becomes invertible. Indeed, the first term in the product  is an involution,
and the second term is invertible by Neumann series arguments.
But then
\[P\xrightarrow{1-P-Q} Q, \quad\text{ and }\quad P\xrightarrow{1-P-Q+2QP} Q.\]

According to this, if $\|\varepsilon\|$ is sufficiently small, then
\[\langle P,\mathrm R\mathrm 0_\Omega\rangle\xrightarrow{\langle1-P-P_\varepsilon+2P_\varepsilon P,1\rangle}
\langle P_\varepsilon,\mathrm R\mathbf 0_\Omega\rangle.\]
After eliminating the unnecessarily stabilization terms in the latter
expression, we see finite dimensionality.
\end{point}
\begin{remark}
When we send $P$ into $Q$ we can choose between $1-P-Q$ and $(1-P-Q)^{-1}$.
If $\mathfrak A$ allows a continuous division by $2$, then
we have more symmetrical choice: We can take the
geometric mean between $1-P-Q$ and $(1-P-Q)^{-1}$. This yields the sign operator
\[\sgn (1-P-Q)=\int \frac{\frac{(1-P)+(1-Q)}2-\frac{P+Q}2z}{\frac{(1-P)+(1-Q)}2+\frac{P+Q}2z}\frac{|dz|}{2\pi}.\]
Here it was understood that
$\frac{P+Q}2= \mathrm R\mathbf 0_\Omega+\frac{(P-\mathrm R\mathrm 0_\Omega)+(Q-\mathrm R\mathbf 0_\Omega)}2$,
and similarly for $\frac{(1-P)+(1-Q)}2$. This yields
\[P\xrightarrow{\sgn (1-P-Q)}Q , \quad\text{ and }\quad
P\xrightarrow{(1-2Q)\sgn (1-P-Q)=\sgn (1-P-Q)\,(1-2P)}Q,\]
which are much more symmetrical choices. This construction fits to the case of $*$-algebras.
\end{remark}
\begin{point}
Using the same techniques, one can prove the following version of Corollary
\ref{cor:tame}: If $\langle b,a\rangle\simalg\langle\tilde b,a\rangle$,
where $b,\tilde b,a$ are finite matrices, $b\approx\tilde b\approx a$,
then there is a finite index set $\Omega'$ and $\psi,\phi\approx1_{\Omega\dot\cup\Omega'}$ such that
$\langle b\oplus\mathrm R\mathbf 0_{\Omega'},a\oplus\mathrm R\mathbf 0_{\Omega'}\rangle
\xrightarrow{\langle\psi,\phi\rangle}
\langle\tilde b\oplus\mathrm R\mathbf 0_{\Omega'},a\oplus\mathrm R\mathbf 0_{\Omega'}\rangle$.
As a consequence, the discussion of individual Grassmannian elments can be reduced to finite matrices.
This is, however, of little help in terms of the global $H$-ring structure.
\end{point}

\section{Fredholm operators}
Here we discuss how Fredholm operators can be placed in this context.
It is instructive to see how standard elements of Fredholm module
theory come up in the context of virtual Grassmannians.
See e.~g.~\cite{GVF} for comparison.
\begin{point} \textbf{Connectors.}
Suppose that $\xi,a\in\Psi_\Omega(\mathfrak A^+)$,
$a$ is an idempotent, and  $\xi$ is an invertible element
such that $\bar a^\xi\approx a$.
In such a pair $[\xi,a\rangle$, we  call the first term as the connector element,
and the second term as the base idempotent.
We define the index
\[\Ind\,[\xi,a\rangle=\langle\bar a^\xi,a\rangle.\]
The nicest case is when $\xi$ is an involution.
We can define the sum $[\xi_1,a_1\rangle\oplus[\xi_2,a_2\rangle=[\xi_1\oplus\xi_2,a_1\oplus a_2\rangle$
and the additive inverse $[\xi,a\rangle^{\inv}=[\xi,\bar a\rangle$.

Suppose now that $c$ is another idempotent,  and $\sigma$ is invertible such that
$\bar c^\sigma\approx c$ but in $\Psi_\Xi(\mathfrak A)$.
Let us consider the product of indices
\[\Ind\,[\xi,a\rangle\overleftarrow\otimes\Ind\,[\sigma,c\rangle=
\langle\bar a^\xi,a\rangle\overleftarrow\otimes\langle\bar c^\sigma,c\rangle=
\langle\bar a^\xi\otimes\bar c^\sigma+a^\xi\otimes c
,a\otimes\bar c^\sigma+\bar a\otimes c \rangle.\]
One can see that
\[a\otimes\bar c+\bar a\otimes c\xrightarrow{a\otimes\sigma+\bar a\otimes 1_\Xi}
a\otimes\bar c^\sigma+\bar a\otimes c\]
and
\[a\otimes c+\bar a\otimes\bar c\xrightarrow{\xi a\otimes 1_\Xi+\xi\bar a\otimes\sigma}
\bar a^\xi\otimes\bar c^\sigma+a^\xi\otimes c.\]
Hence, it is reasonable  to  define
\[[\xi,a\rangle\overleftarrow\otimes[\sigma,c\rangle:=
\bigl[(a\otimes\sigma+\bar a\otimes 1_\Xi)^{-1}(\xi a\otimes 1_\Xi+\xi\bar a\otimes\sigma),
a\otimes\bar c+\bar a\otimes c\bigr\rangle.\]
Then the product of indices is equivalent to the index of the product. Indeed:
\[\Ind\,[\xi,a\rangle\overleftarrow\otimes[\sigma,b\rangle
\xrightarrow{\langle a\otimes\sigma+\bar a\otimes 1_\Xi,a\otimes\sigma+\bar a\otimes 1_\Xi\rangle}
\Ind\,[\xi,a\rangle\overleftarrow\otimes\Ind\,[\sigma,b\rangle.\]
A bonus is that the connector term of $[\xi,a\rangle\overleftarrow\otimes[\sigma,c\rangle$
is an involution if $\xi,\sigma$ are involutions.
\end{point}
\begin{lemma} The following hold:

(a) $[\xi_1,a_1\rangle\approx[\xi_2,a_2\rangle\Rightarrow\Ind\,[\xi_1,a_1\rangle\simalg\Ind\,[\xi_2,a_2\rangle$;

(b) $\xi\bar a\xi^{-1}=a\Rightarrow\Ind\,[\xi,a\rangle \simalg \mathbf 0$;

(c) $\Ind\,[\xi,a\rangle \simalg \Ind\,[\xi^{-1},a\rangle$,
$\langle\xi^2 a \xi^{-2},a \rangle \simalg \mathbf 0$;

(d) $\Ind\,[\xi_1,a_1\rangle\oplus[\xi_2,a_2\rangle\simeq\Ind\,[\xi_1,a_1\rangle\oplus\Ind\,[\xi_2,a_2\rangle$;

(e) $\Ind\,[\xi,a\rangle^{\inv}\simeq \left(\Ind\,[\xi,a\rangle\right)^{\inv}$;

(f) $\Ind\,[\xi,a\rangle\overleftarrow\otimes[\sigma,b\rangle\simalg
\Ind\,[\xi,a\rangle\overleftarrow\otimes\Ind\,[\sigma,b\rangle$.
\begin{proof} The points (b), (d), (e), (f) are immediate.

(a${_1}$) $\xi_1\approx\xi_2\Rightarrow\Ind\,[\xi_1,a\rangle\simalg \Ind\,[\xi_2,a\rangle$
follows from $\Ind\,[\xi_1,a\rangle\xrightarrow{\langle\xi_2\xi_1^{-1},1_\Omega\rangle}\Ind\,[\xi_2,a\rangle$.

(a${_2}$) $a_1\approx a_2\Rightarrow\Ind\,[\xi,a_1\rangle\simalg \Ind\,[\xi,a_2\rangle$ comes as follows:
The algebraic homotopies
\[\Ind\,[\xi,a_1\rangle\oplus[\xi,\bar a_2\rangle
\xrightarrow{\langle\sw_{12}(a_2), \sw_{12}(a_2)\rangle}
\Ind\,[ \sw_{12}(a_2)(\xi\oplus\xi)\sw_{12}(a_2),
\mathrm R\langle a_1,a_2\rangle\rangle\simalg \]
(by (a${}_1$))
\[\simalg \Ind\,\left[\begin{bmatrix}&\xi\\\xi&\end{bmatrix},\mathrm R\langle a_1,a_2\rangle\right\rangle
=\langle \mathrm R\langle a_1^\xi,a_2^\xi\rangle,\mathrm R\langle a_1,a_2\rangle\rangle
\simalg \]
(by Corollary \ref{cor:diff})
\[\simalg \mathbf R\langle a_1^\xi,a_2^\xi\rangle\oplus\mathbf R\langle a_1,a_2\rangle^{\inv}
\simalg \mathbf R\langle a_1,a_2\rangle\oplus\mathbf R\langle a_1,a_2\rangle^{\inv}\simalg \mathbf 0\]
imply the statement.

(c) According to (a${_2}$), we see that
$\Ind\,[\xi,a\rangle\simalg \Ind\,[\xi,\xi^{-1}\bar a\xi\rangle
\simalg \Ind\,[\xi,\xi^{-1}\bar a\xi\rangle'=\Ind\langle\xi^{-1},a\rangle.$
But then, according to Corollary \ref{cor:diff}, we also see that
$\langle\xi a\xi^{-1},\xi^{-1} a\xi\rangle\simalg \mathbf 0$. Consequently,
$\langle\xi^2 a\xi^{-2},a\rangle\simalg \mathbf 0$.
\end{proof}
\end{lemma}
\begin{point}\textbf{Fredholm connectors.}
In general, $\xi$ can always be replaced by an involution.
Indeed, suppose  that $a$ is an idempotent, $\xi a=\bar a\xi$, and $\xi$ allows
a parametrix $\eta$ such that $\xi\eta\approx\eta\xi\approx1$.
Then it is not hard to see that, for example,
\[\tilde\xi=(a+a\xi\bar a-\bar a)(a-a\eta\bar a-\bar a)(a+a\xi\bar a-\bar a)\]
will be an involution. This is a reasonable substitute.
\begin{lemma}
If $\xi$ is invertible, then $\Ind\,[\xi,a\rangle\simalg\Ind\,[\tilde\xi,a\rangle$.
\begin{proof}
Apply the following lemma with $\psi=\tilde\xi\xi^{-1}$.
\end{proof}
\end{lemma}
\begin{lemma}
Suppose that $\psi$ is invertible such that $\psi-\bar a\psi \bar a\approx a$
 or  $\psi-a\psi a\approx\bar a$. Then
\[\langle \psi a\psi ^{-1} ,a\rangle\simalg\mathbf 0.\]
\begin{proof}
Consider the first case. Then
\[\langle \psi a\psi ^{-1} ,a\rangle\simalg
\left\langle \begin{bmatrix}\psi&\\&1 \end{bmatrix}\begin{bmatrix}a&\\&0 \end{bmatrix}
\begin{bmatrix}\psi^{-1}&\\&1 \end{bmatrix}
,\begin{bmatrix}a&\\&0 \end{bmatrix}\right\rangle
\rightarrow\left\langle\begin{bmatrix}a&\\&0 \end{bmatrix},\begin{bmatrix}a&\\&0 \end{bmatrix}\right\rangle
\simalg\mathbf 0,\]
where the conjugation in the middle is induced by the replacement in the conjugating term according to the
observation
\[\begin{bmatrix}\psi&\\&1 \end{bmatrix}\approx
\begin{bmatrix}\bar a\psi\bar a+ a&\bar a\psi a\\a\psi\bar a&\bar a+a\psi a \end{bmatrix}
=\begin{bmatrix}\bar a&a\\a&\bar a \end{bmatrix}
\begin{bmatrix}\psi&\\&1 \end{bmatrix}
\begin{bmatrix}\bar a&a\\a&\bar a \end{bmatrix}.\]
The second case is similar.
\end{proof}
\end{lemma}
The price we paid for having $\tilde\xi$  is that we have departed from the case of
conjugating matrices which are manifestly unitary in the $*$-algebraic case.
However, it is this setting which is applied in the discussion of the index of Fredholm operators as follows.
\end{point}
\begin{point} \textbf{Fredholm operators.}
We say that $\psi\in\Psi_{\Omega_1,\Omega_0}(\mathfrak A^+)$ is a Fredholm operator if
there is a parametrix  $\phi\in\Psi_{\Omega_0,\Omega_1}(\mathfrak A^+)$ such that
$\phi\psi\approx 1_{\Omega_0}$ and $\psi\phi\approx 1_{\Omega_1}$.
We say that $\xi\in \Psi_{\Omega_0\cup\Omega_1}(\mathfrak A^+)$ is a connector element
for the Fredholm operator $\psi$, if it is invertible and
\[\xi\approx\begin{bmatrix}&\phi\\\psi&\end{bmatrix}.\]
One can easily obtain such an element by taking the involution
\[\mathrm F(\psi,\phi)
=\begin{bmatrix}1_{\Omega_0}-\phi\psi&2\phi-\phi\psi\phi\\\psi&\psi\phi-1_{\Omega_1}\end{bmatrix}=
\begin{bmatrix}1_{\Omega_0}&\phi\\&-1_{\Omega_1}\end{bmatrix}
\begin{bmatrix}1_{\Omega_0}&\\-\psi&-1_{\Omega_1}\end{bmatrix}
\begin{bmatrix}1_{\Omega_0}&\phi\\&-1_{\Omega_1}\end{bmatrix}.\]
\end{point}
\begin{point}
To a connector element $\xi$ as above, we can associate the index
\[\Ind\,[\xi,0_{\Omega_0}\oplus1_{\Omega_1}\rangle
=\left\langle\xi\begin{bmatrix}1_{\Omega_0}&\\&0_{\Omega_1}\end{bmatrix}\xi^{-1},
\begin{bmatrix}0_{\Omega_0}&\\&1_{\Omega_1}\end{bmatrix}\right\rangle.\]
It is easy to see that, up to smooth conjugation, this Grassmannian element depends only on $\psi$.
Nevertheless, in order to get a well-defined index element,
we define the index of the Fredholm pair $(\psi,\phi)$ as
\[\Ind_{\mathrm F} (\psi,\phi):=\Ind\,[\mathrm F(\psi,\phi),0_{\Omega_0}\oplus 1_{\Omega_1}\rangle.\]
\end{point}
\begin{point}
We can define the sum of Fredholm pairs as
$(\psi_1,\phi_1)\oplus(\psi_2,\phi_2):=(\psi_1\oplus\psi_2,\phi_1\oplus\phi_2)$,
and the inverse pair as
$(\psi,\phi)^{\inv}=(\phi,\psi)$.
If $(\psi',\phi')$ is a Fredholm pair such that
$\psi'\in\Psi_{\Omega_2,\Omega_1}, \phi'\in\Psi_{\Omega_1,\Omega_2}$,
then we define the composition  as
$(\psi',\phi')\circ(\psi,\phi)=(\psi'\psi,\phi\phi').$
\end{point}
\begin{lemma}\label{lem:Fredholm} The following hold:

(a) $(\psi_1,\phi_1)\approx(\psi_2,\phi_1)\Rightarrow
\Ind_{\mathrm F} (\psi_1,\phi_1)\simalg \Ind_{\mathrm F} (\psi_2,\phi_2)$;

(b)  if $\psi$ or $\phi$ is invertible, then
$ \Ind_{\mathrm F} (\psi,\phi) \simalg \mathbf 0$;

(c) $\Ind_{\mathrm F} (\psi_1,\phi_1)\oplus(\psi_2,\phi_2)\simeq \Ind (\psi_1,\phi_1)\oplus\Ind (\psi_2,\psi_2)$;

(d) $\Ind_{\mathrm F}  (\psi,\phi)^{\inv}\simeq \left(\Ind_{\mathrm F}  (\psi,\phi)\right)^{\inv}$;

(e) $\Ind_{\mathrm F}  (\psi',\phi')\circ(\psi,\phi)\simalg
\Ind_{\mathrm F}  (\psi,\phi)\oplus\Ind_{\mathrm F}  (\psi',\psi')$;

(f) $\Ind_{\mathrm F} (\psi,\phi)\overleftarrow\otimes(\theta,\chi)\simalg
\Ind_{\mathrm F} (\psi,\phi)\overleftarrow\otimes\Ind_{\mathrm F} (\theta,\chi)$.
\begin{proof}[Proof, except of (f)]
Everything is straightforward, except (e).
It is easy to check if $(\psi',\phi')=(u, u^{-1})$, where $u$ is invertible.
In this case
\[\Ind_{\mathrm F} (\psi,\phi)\xrightarrow{\langle 1_{\Omega_0}\oplus u,1_{\Omega_0}\oplus u\rangle}
\Ind_{\mathrm F}  (\psi',\phi')\circ(\psi,\phi).\]
This special case implies the general case as follows.
Let $\xi'=\begin{bmatrix}\xi'_{11}&\xi'_{12}\\\xi'_{21}&\xi'_{22}\end{bmatrix}
\approx\begin{bmatrix}&\phi'\\\psi'&\end{bmatrix}$
be a connector element with inverse $\tilde\xi'$. Then
\[\left(\begin{bmatrix}\xi'_{21}&\xi'_{22}\\\xi'_{11}&\xi'_{12}\end{bmatrix},
\begin{bmatrix}\tilde\xi'_{12}&\tilde\xi'_{11}\\\tilde\xi'_{22}&\tilde\xi'_{21}\end{bmatrix}\right)
\circ\left(\begin{bmatrix}\psi&\\&\psi'\end{bmatrix},\begin{bmatrix}\phi&\\&\phi'\end{bmatrix}\right)
\approx\left(\begin{bmatrix}\psi'\psi&\\&1_{\Omega_1}\end{bmatrix},
\begin{bmatrix}\phi'\phi&\\&1_{\Omega_1}\end{bmatrix}\right).\]
Consequently,
$\Ind_{\mathrm F}  (\psi,\phi)\oplus\Ind_{\mathrm F}  (\psi',\psi')\simalg
\Ind_{\mathrm F}  (\psi',\phi')\circ(\psi,\phi)\oplus \Ind_{\mathrm F}  (1_{\Omega_1},1_{\Omega_1}).$
\end{proof}
\end{lemma}

What prevents us from proving (f) is that we have not defined the product of Fredholm
pairs yet. Thus, it remains to make a definition such that (f) would hold.

\begin{point}Evaluating it, we find that the connector term of
$[\mathrm F(\psi,\phi),0_{\Omega_0}\oplus 1_{\Omega_1}\rangle\overleftarrow\otimes
[\mathrm F(\theta,\chi),0_{\Xi_0}\oplus 1_{\Xi_1}\rangle$
is a connector term for the Fredholm pair
\[\left(\begin{bmatrix}(1_{\Omega_0}-\phi\psi)\otimes\theta&(2\phi-\phi\psi\phi)\otimes1_{\Xi_1}\\
\psi\otimes1_{\Xi_0}&(\psi\phi-1_{\Omega_1})\otimes\chi\end{bmatrix},
\begin{bmatrix}(1_{\Omega_0}-\phi\psi)\otimes\chi&(2\phi-\phi\psi\phi)\otimes1_{\Xi_0}\\
\psi\otimes1_{\Xi_1}&(\psi\phi-1_{\Omega_1})\otimes\theta\\\end{bmatrix}\right),\]
which, henceforth, can be chosen as the definition of $(\psi,\phi)\overleftarrow\otimes(\theta,\chi)$.
Applying the composition
\[\left(\begin{bmatrix}1_{\Omega_0\times\Xi_1}&-\phi\otimes\theta\\&-1_{\Omega_1\times\Xi_0}\end{bmatrix},
\begin{bmatrix}1_{\Omega_0\times\Xi_1}&-\phi\otimes\theta\\&-1_{\Omega_1\times\Xi_0}\end{bmatrix}\right)\circ,\]
which, according to \ref{lem:Fredholm} (b) and (e), does not change the index
up to algebraic homotopy, we obtain $\approx$
\[\left(\begin{bmatrix}1_{\Omega_0}\otimes\theta&\phi\otimes1_{\Xi_1}\\
-\psi\otimes1_{\Xi_0}&(1_{\Omega_1}-\psi\phi)\otimes \chi\end{bmatrix},
\begin{bmatrix}(1_{\Omega_0}-\phi\psi)\otimes\chi&-\phi\otimes 1_{\Xi_0}\\
\psi\otimes1_{\Xi_1}&1_{\Omega_1}\otimes\theta\end{bmatrix}\right). \]

Hence, this is another sufficiently good definition for $(\psi,\phi)\overleftarrow\otimes(\theta,\chi)$.
This latter one is essentially the same as \cite{F} (2.6). Needless to say, many variants are possible.
\end{point}
Fredholm operators coming from geometry are more general.
They do not act around $0_{\Omega_0}\oplus 1_{\Omega_1}$
but more general idempotents corresponding to vector bundles.

\end{document}